\documentclass{amsart}
\pdfoutput=1
\usepackage[dvipsnames]{xcolor}
\usepackage[unicode,
  colorlinks=true,
  linktocpage=true,
  citecolor=Emerald,
  linkcolor=Fuchsia,
  urlcolor=Fuchsia,
  pdftitle={Partial groups as symmetric simplicial sets},
  pdfauthor={Philip Hackney and Justin Lynd},
  pdfkeywords={partial group, groupoid, symmetric simplicial set}]{hyperref}
\usepackage{tikz-cd}
\usepackage[capitalise]{cleveref}
\usepackage{enumitem}
\usepackage{mathtools}
\usepackage{mathrsfs}
\usepackage{amssymb}
\usepackage{microtype}
\newcommand{\sset}{\mathsf{sSet}}
\newcommand{\symset}{\mathsf{Sym}}
\newcommand{\set}{\mathsf{Set}}
\newcommand{\op}{{\operatorname{op}}}
\DeclareMathOperator{\id}{id}
\newcommand{\fplus}{\mathbf{\Upsilon}}
\newcommand{\simpcat}{\mathbf{\Delta}}
\newcommand{\rep}[1]{\Upsilon^{#1}}
\newcommand{\finplus}{\mathsf{Fin}_+}

\DeclareMathOperator*{\colim}{colim}

\newcommand{\edgemap}{\mathscr{E}}
\DeclareMathOperator{\mat}{Mat}
\newcommand{\gpd}{\mathsf{Gpd}}
\newcommand{\pgpd}{\mathsf{pGpd}}
\newcommand{\pgrp}{\mathsf{pGrp}}
\newcommand{\fpg}[1]{\mathfrak{F}^{#1}}

\newcommand{\spine}[1]{\textup{Sp}^{#1}}

\numberwithin{equation}{section}
\newtheorem{theorem}{Theorem}[section]
\newtheorem{lemma}[theorem]{Lemma}
\newtheorem{proposition}[theorem]{Proposition}
\newtheorem{corollary}[theorem]{Corollary}
\theoremstyle{definition}
\newtheorem{definition}[theorem]{Definition}
\newtheorem{example}[theorem]{Example}
\newtheorem{remark}[theorem]{Remark}
\newtheorem{notation}[theorem]{Notation}
\theoremstyle{remark}
\newtheorem{warning}[theorem]{Warning}

\title{Partial groups as symmetric simplicial sets}

\author{Philip Hackney}
\address{Department of Mathematics, University of Louisiana at Lafayette, USA}
\email{philip@phck.net} 
\urladdr{http://phck.net}

\author{Justin Lynd}
\address{Department of Mathematics, University of Louisiana at Lafayette, USA}
\email{lynd@louisiana.edu}

\thanks{This work was supported by a grant from the Simons Foundation (\#850849, PH). JL was partially supported by NSF Grant DMS-1902152.}

\keywords{partial group, groupoid, symmetric simplicial set}

\subjclass[2020]{Primary: 
18N50, % Simplicial sets, simplicial objects
20N99, % Other generalizations of groups
18A30, % Limits and colimits
18A40. % Adjoint functors 
Secondary: 
18F20, % Presheaves and sheaves, stacks, descent conditions
20D20, % Sylow subgroups, Sylow properties, π-groups, π-structure
20L05, % Groupoids
55U10} % Simplicial sets and complexes in algebraic topology

\begin{document}

\begin{abstract}
We give a new characterization of partial groups as a  subcategory of symmetric (simplicial) sets.
This subcategory has an explicit reflection, which permits one to compute colimits in the category of partial groups.
We also introduce the notion of a partial groupoid, which encompasses both groupoids and partial groups.
\end{abstract}

\maketitle

\setcounter{tocdepth}{1}
\tableofcontents

A symmetric set is a simplicial set equipped with the additional action of symmetric groups at each level.
Grothendieck showed that the category of groupoids is a full subcategory of symmetric sets, whose essential image consists of those $X$ for which the Grothendieck--Segal map
\[ \begin{tikzcd}[sep=small]
X_n \rar{\cong} & X_1 \times_{X_0} \cdots \times_{X_0} X_1  
\end{tikzcd} \]
is a bijection for every $n\geq 1$ \cite[4.1]{Grothendieck:TCTGA3}.
Elements of $X_n$ should be interpreted as length $n$ strings of composable morphisms in the associated groupoid.
Likewise, the category of groups may be viewed as a full subcategory of the reduced symmetric sets, which are those where $X_0$ has a single element.

Partial groups were introduced in \cite{Chermak:FSL} for the purpose of
constructing group-theoretic avatars of transporter systems, which are finite
categories modeling the $p$-local structure of a finite group.  
In a group $G$, one can form the total multiplication of any finite word in the underlying set $G$. 
In a partial group, this is weakened so that only certain admissible words admit a total multiplication; a group is just a partial group where every word is admissible. 
Partial groups were recently shown by Lemoine and Molinier \cite{LemoineMolinier:PGPRFS} to subsume the pregroups of Stallings \cite{Stallings:GT3DM}.
In \cite{gonzalez}, Gonz\'alez showed that partial groups may be considered as simplicial sets satisfying certain properties (see \cref{gonzalez theorem} below), the perspective that we will adopt in this paper.
Gonz\'alez's characterization of partial groups extends in the expected way to give a notion of \emph{partial groupoids}.\footnote{Elsewhere in the literature, the term \emph{partial groupoid} has a different meaning: a set together with a partially-defined binary operation.}
As with the passage from groups to groupoids, this added flexibility is very useful and simplifies the formal presentation below.

A \emph{spiny} symmetric set (\cref{def spiny symset}) is one in which each $n$-simplex is uniquely determined by an appropriately chosen $n$-tuple of its edges.
{
\renewcommand{\thetheorem}{A}
\begin{theorem}
The category of spiny symmetric sets is equivalent to the category of partial groupoids.
The category of reduced spiny symmetric sets is equivalent to the category of partial groups.
\end{theorem}
}
This characterization (\cref{thm spiny equals pgpd} and \cref{cor spiny partial groups} below) allows for a very concise definition of partial groupoid or partial group, as a known type of object satisfying one or two simple properties.
It also puts inversion in its proper place: in the characterization of partial groups as simplicial sets, the inversion is specified as additional structure of certain simplicial sets, when in reality it is a property they may possess. 

In 1988, Lawvere observed that the category of groupoids is a reflective subcategory of the category of symmetric sets \cite{Lawvere:TGCDCTFA}.
This means that the inclusion functor admits a left adjoint, called the reflection.
We extend this result and show that the larger category of partial groupoids forms a reflective subcategory of symmetric sets.
Unlike the situation for the reflection to groupoids, no composite morphisms must be added, so the spiny reflection of a symmetric set is a quotient. 
{
\renewcommand{\thetheorem}{B}
\begin{theorem}\label{reflection theorem}
The category of (reduced) spiny symmetric sets is a reflective subcategory of the category of symmetric sets.
\end{theorem}
}

As a consequence, we recover the following known theorem, whose two parts are due to Chermak and Salati, respectively.
Existence of limits first appeared in \cite[Appendix A]{Chermak:FL}, while existence of colimits is the primary subject of \cite{Salati:LCGRPG}.
{
\renewcommand{\thetheorem}{C}
\begin{corollary}\label{cor bicomplete}
The category of partial groups is complete and cocomplete.  
\end{corollary}
}
Likewise, the category of partial groupoids is bicomplete.
Limits and colimits of symmetric sets are constructed levelwise, as are limits of (reduced) spiny symmetric sets.
Colimits of (reduced) spiny symmetric sets are more delicate, but are constructed by a standard process: first compute the colimit in the larger category of symmetric sets, and then apply the reflection.

Throughout the paper we include several examples of partial groups and colimits of such. 
\Cref{ex nilpotence commutative} gives partial groups related to commutativity and nilpotence in a group $G$.
\Cref{ex reduction groupoid} explains how every groupoid gives rise to a partial group.

A theoretically important collection of partial groups are the \emph{word classifiers} $\fpg{m}$ which we construct (as symmetric sets) in \cref{example reduced representable}.
These are partial groups with the property that for any partial group $X$, the set $\hom(\fpg{m}, X)$ is precisely the set of length $m$ admissible words in $X$.
Every partial group is canonically a colimit of word classifiers. 
To our knowledge these partial groups have not explicitly appeared in the literature before, except for $m = 1$, though they are simple to construct in our framework. 

In \cref{ex colimits} we discuss coproducts of partial groups and partial groupoids.
In contrast to the coproduct (free product) of groups, coproducts of partial groups are quite simple: they are wedge sums. 
Pushouts of partial groups are occasionally simple as well. 
As we explain in \cref{ex pushout along G}, the amalgamated union of two partial groups sharing a common subgroup is easy to construct, as it agrees with the pushout in symmetric sets.
It is important that the intersection of the partial groups is a group and not a more general partial group: in \cref{ex counterexample} we give an example where the pushout does not agree with the one in symmetric sets.

\section{Simplicial sets, partial groups, and partial groupoids}

\subsection{Simplicial sets and symmetric sets}
The simplicial category $\simpcat$ has objects the sets $[n] = \{ 0, 1, \dots, n \}$ for $n\geq 0$, and morphisms are the (weakly) order-preserving maps. 
A \emph{simplicial set} is a presheaf on $\simpcat$, that is, a functor $X \colon \simpcat^\op \to \set$.
We write $\sset$ for the category of simplicial sets and their morphisms (natural transformations of functors).
As is usual, we denote by $X_n$ the value of $X$ at the object $[n]$. 
There are generating morphisms for $\simpcat^\op$ denoted, for $0\leq i \leq n$, by $d_i \colon [n] \to [n-1]$ (face operators) and $s_i \colon [n] \to [n+1]$ (degeneracy operators) satisfying certain relations, and often a simplicial set is written as
\[ \begin{tikzcd}
X_0 \rar["s_0" description] & X_1 \lar[shift left=2, "d_1"] \lar[shift right=2, "d_0"']  \rar[shift left=1.5] \rar[shift right=1.5] & X_2 \lar[shift left=3,"d_2"] \lar[shift right=3,"d_0"'] \lar
\rar[shift left=3] \rar[shift right=3] \rar &
X_3 
\lar[shift left=1.5] \lar[shift left=4.5] \lar[shift right=1.5] \lar[shift right=4.5] \cdots
\end{tikzcd} \]
prioritizing these generating maps.
The map $d_i \colon X_n \to X_{n-1}$ is the value of the unique injection $[n-1] \to [n]$ missing the element $i$, and $s_i \colon X_n \to X_{n+1}$ is the value of the unique surjection $[n+1] \to [n]$ hitting the element $i$ twice.
As usual, if $\alpha \colon [n] \to [m]$ is an arbitrary order-preserving map, we tend to write $\alpha^* \colon X_m \to X_n$ instead of $X(\alpha)$ for the value of $X$ at $\alpha$.

\begin{example}\label{ex nerve of category}
If $C$ is a (small) category, then its nerve, denoted $NC$, is a simplicial set.
The set of vertices $NC_0$ is the set of objects of $C$, while the set of edges $NC_1$ is the set of morphisms of $C$.
General $n$-simplices are chains of $n$ composable morphisms
\[ \begin{tikzcd}[sep=small]
x_0 \rar{f_1} & x_1 \rar{f_2} & x_2 \rar{f_3} & \cdots \rar{f_{n-1}} & x_{n-1} \rar{f_n} & x_n.
\end{tikzcd} \]
Simplicial sets isomorphic to nerves of categories can be described by a certain exactness condition, called the \emph{Grothendieck--Segal condition}. (See \S\ref{simp sets pgpd}.)
\end{example}

The category $\fplus$ has the same objects as $\simpcat$, but morphisms in $\fplus$ are arbitrary functions between these sets.\footnote{Unlike $\simpcat$, there is no standard name for this category: we have seen more than ten names in the literature, including $\mathbb{S}$, $\mathsf{F}$, $\mathbf{F}$, $\widetilde{\mathbf{\Delta}}$, $!\mathbf{\Delta}$, $\mathbf{\Delta}_{\mathrm{sym}}$, and others. Our choice follows \cite{Cisinski:PCMTH}.} 
A \emph{symmetric (simplicial) set} is a presheaf on $\fplus$, and we write $\symset = [\fplus^\op, \set]$ for the category of symmetric sets.
Symmetric sets are simplicial sets with additional structure, namely compatible actions of the symmetric groups $\Sigma_{n+1}^\op$ on $X_n$ for each $n$, see \cite[\S3]{Grandis:FSSSS} for a description via generators and relations.
We provide an alternative set of generating morphisms in \cref{prop fplus generation}.
Symmetric sets are nearly always called symmetric \emph{simplicial} sets in the literature, but to avoid excessive alliteration later we always use this abbreviated name.\footnote{Readers are warned to not carry this too far: one would not want to call a functor $\fplus^\op \to \mathsf{Grp}$ a \emph{symmetric group}.}
As $\fplus$ is equivalent to the category $\finplus$ of non-empty finite sets, occasionally we will regard symmetric sets as functors $\finplus^\op \to \set$.
There is an inclusion functor $J\colon \simpcat \to \fplus$, which gives a (faithful) restriction $J^* \colon \symset \to \sset$ forgetting about the maps $\phi^* \colon X_m \to X_n$ whenever $\phi \colon [n] \to [m]$ is not order preserving (in particular, the actions by symmetric groups are forgotten).

The category of symmetric sets is very rich; for instance it admits a Quillen model structure for $\infty$-groupoids, see \cite[8.3.8]{Cisinski:PCMTH}.
As mentioned in the introduction, 
Grothendieck gave a nerve functor for groupoids, landing in the category of symmetric sets \cite[4.1]{Grothendieck:TCTGA3}.\footnote{Terminology has changed over time, and we warn that in \cite{Grothendieck:TCTGA3} \emph{ensemble simplicial} means \emph{symmetric set} and \emph{ensemble semi-simplicial} means \emph{simplicial set}.}
See \cite[\S4]{BergerMelliesWeber:MAAT} for more details about the relationship between symmetric sets and groupoids, and \cite[\S2]{Grandis:SHFFSS} for an explanation of the reflection $\symset \to \gpd$.
If $C$ is a (small) category, then its nerve $NC\in \sset$ is the restriction of a symmetric set if and only if $C$ is a groupoid, that is, if every morphism is invertible.
\Cref{rmk cyclic} is a slightly more refined result, which we will not need in what follows.

\begin{warning}
There is a related, but different notion in the literature, namely that of a presheaf over the larger category $\simpcat \mathbf{S}$. 
A closely related notion to this one is that of a cyclic set in the sense of Connes \cite{Connes:CCFE}, namely a presheaf on the category $\simpcat \mathbf{C}$ or $\mathbf{\Lambda}$ (see \cite{Loday:CH} for both of these).  
There are functors $\simpcat \to \mathbf{\Lambda} \to \simpcat \mathbf{S} \to \fplus$ which are the identity on objects. 
The first two functors are faithful, whereas the last one is full but not faithful.
Nevertheless, the composite is the usual inclusion.
\end{warning}

\begin{remark}[Groupoids and cyclic sets]\label{rmk cyclic}
A category $C$ is a groupoid if and only if its nerve is the restriction of a symmetric set in either sense.
In fact, if $NC$ is the restriction of a cyclic set,
then $C$ is a groupoid (the converse holds using the previous statement and the
factorization $\simpcat \to \mathbf{\Lambda} \to \fplus$).  One can show that an
arbitrary morphism  $f \colon a \to b$ is invertible.  Calculating $d_0t_2$ and
$d_1t_2$ applied to the 2-simplex
\[ \begin{tikzcd}
a \rar{f} & b \rar{t_1(\id_b)} & b
\end{tikzcd} \]
using the cyclic identities (see \cite[\S6.1]{Loday:CH}) establishes that $f$ has a right inverse.
Calculating $d_2t_2t_2$ and $d_1t_2t_2$ applied to the 2-simplex 
\[ \begin{tikzcd}
a \rar{t_1(\id_a)} & a \rar{f} & b
\end{tikzcd} \] 
establishes that $f$ has a left inverse.
Since every morphism is invertible, $C$ is a groupoid.
(Be warned that $NC$ could admit more than one cyclic structure.)
\end{remark}
 
\begin{notation}
Certain notation follows \cref{ex nerve of category}.
If $X$ is a simplicial set and $x \in X_0$, we write $\id_x \in X_1$ for $s_0(x)$.
An edge $f\in X_1$ has a source $d_1f \in X_0$ and a target $d_0f \in X_0$, and we write $f \colon x_0 \to x_1$ where $x_0 = d_1f$ and $x_1 = d_0f$.
\end{notation}

A \emph{reduced} simplicial set is a simplicial set where $X_0$ is a point, and likewise for a reduced symmetric set.

\begin{notation}[Special maps]
We name several important maps in $\fplus$.
\begin{itemize}
\item The flip maps $\tau_n \colon [n] \to [n]$ which act by $i \mapsto n-i$.
\item The folding maps $\chi_n \colon [2n] \to [n]$ which act by $i \mapsto |n-i|$.
\item 
The edge classifiers $\rho_{ij} = \rho_{i,j} \colon [1] \to [n]$ which are given by $\rho_{ij}(0) = i$ and $\rho_{ij}(1) = j$. Note that $\rho_{ij} \tau_1 = \rho_{ji}$.
\end{itemize}
The map $\rho_{ij}$ is in $\simpcat$ if and only if $i\leq j$, and $\chi_n$, $\tau_n$ are in $\simpcat$ if and only if $n=0$ (in which case they are identities).
\end{notation}

\begin{definition}\label{def opposite}
If $X$ is a simplicial set, then $X^\op \in \sset$ denotes its opposite. 
It has the same underlying sequence of sets.
If $\alpha \colon [n] \to [m]$ is an order-preserving map, then so is $\tau_m \alpha \tau_n \colon [n] \to [m]$, and $X^\op(\alpha) \colon X_m \to X_n$ is defined to be  $X(\tau_m \alpha \tau_n)$.
\end{definition}
In particular, $X^\op(\rho_{ij}) = X(\rho_{n-j,n-i}) \colon X_n \to X_1$, and $(NC)^\op = N(C^\op)$.

\subsection{Simplicial sets and partial group(oid)s}\label{simp sets pgpd}

\begin{definition}\label{def edgy sset}
A simplicial set $X$ is called \emph{edgy} if
\[
	\edgemap = \edgemap_n = \prod_{i=1}^n \rho_{i-1,i}^* \colon X_n \to \prod_{i=1}^n X_1
\]
is injective for $n\geq 1$.
\end{definition}
If $C$ is a (small) category, then its nerve $NC$ is edgy:
\[ \begin{tikzcd}[sep=small]
NC_n \rar{\cong} & NC_1 \times_{NC_0} \cdots \times_{NC_0} NC_1 \rar[hook] & \displaystyle\prod_{i=1}^n NC_1
\end{tikzcd} \]
The first map in this display is an isomorphism, while the second is an injection.
Notice that if $X$ is edgy then the Grothendieck--Segal map $X_n \to X_1 \times_{X_0} \cdots \times_{X_0} X_1$ is an injection, but is not necessarily a bijection.
Because of their relationship to nerves of categories, edgy simplicial sets were called N-simplicial sets in \cite{gonzalez}.

We often adopt a notation like
\begin{equation}
[f_1 | f_2 | \cdots | f_n]
\label{standard form}
\end{equation}
for an $n$ simplex of an edgy simplicial set $X$ when $n>0$.
Here, $f_i$ is the image in $X_1$ under $\rho_{i-1,i}^*$, so this string uniquely determines the simplex.
As we have  $d_0(f_{i}) = d_1(f_{i+1})$ and could instead write this $n$-simplex as 
\begin{equation}\label{eq expanded simplex} \begin{tikzcd}
x_0 \rar{f_{1}} & x_1 \rar{f_{2}} & \cdots \rar{f_{n-1}} & x_{n-1} \rar{f_n} & x_n,
\end{tikzcd} \end{equation}
where the $x_i \in X_0$ are objects.

\begin{definition}
Suppose $X$ is a simplicial set. 
An \emph{anti-involution} on $X$ is a map
\[
	\nu \colon X^\op \to X
\]
such that $\nu^\op \circ \nu = \id_{X^\op}$.
\end{definition}

In particular, we have involutions $\nu_n \colon X_n \to X_n$ for all $n$. 
We write $f\mapsto f^\dagger$ for the involution on $X_1$.
If $X$ is an edgy simplicial set equipped with an anti-involution, then $\nu_n \colon X_n \to X_n$ takes the form
\[
	[f_1 | \cdots | f_n] \mapsto [f_n^\dagger | \cdots | f_1^\dagger].
\]

\begin{definition}[Partial Groupoids]\label{def partial groupoids}
Let $X$ be an edgy simplicial set equipped with an anti-involution $\nu$ such that $\nu_0 \colon X_0 \to X_0$ is the identity map.
We say that $X$ is a \emph{partial groupoid} if for each $n\geq 1$ there is a function $L \colon X_n \to X_{2n}$ such that
\[
	L[f_1 | \cdots | f_n] = [f_n^\dagger | \cdots | f_1^\dagger | f_1 | \cdots | f_n]
\]
and the diagram below-left commutes (whose bottom portion is associated to the order-preserving maps below-right).
\[ \begin{tikzcd}
X_n \rar{L} \dar[swap]{d_0^n} & X_{2n} \dar{\rho_{0,2n}^*} & {[n]} & {[2n]} \\
X_0 \rar[swap]{s_0} & X_{1} & {[0]} \uar{n} & {[1]} \lar \uar[swap]{0,2n}
\end{tikzcd} \]
In other words, the total composite of $L[f_1 | \cdots | f_n]$ is the identity on $x_n = d_0(f_n)$.
\end{definition}

The non-reduced feature of this definition appears to be new. 
The designation `partial groupoid' is in reference to the following, which appears as Theorem 4.8 of \cite{gonzalez}.
It turns out that maps of simplicial sets between partial groupoids automatically preserve the anti-involution (see \cref{lem inv compat automatic}); this is already implicit in Gonz\'alez's theorem.

\begin{theorem}[Gonz\'alez]\label{gonzalez theorem}
The full subcategory of simplicial sets on the reduced partial groupoids is equivalent to the category of partial groups.
\end{theorem}

\begin{example}\label{ex nilpotence commutative}
Suppose $G$ is a (discrete) group.
We temporarily write $BG$ for the usual simplicial set with $BG_n = G^{\times n}$; in this way, we can regard every group as a partial group. 
But we could instead consider the simplicial subset $B_{\textup{com}}G$, with $(B_{\textup{com}}G)_n = \hom(\mathbb{Z}^n,G) \subseteq G^{\times n}$ the subset of pairwise commuting elements (see \cite[\S2]{AdemGomez:CSCLG}).
Then the simplicial set $B_{\textup{com}}G$ is a partial group, which is only a `group' if $G$ is abelian.
Two elements $g,h \in (B_{\textup{com}}G)_1 = G$ can be multiplied only if $gh=hg$.
More generally, let $\mathcal{F}_n = \langle r_1, \dots, r_n \rangle$ be the free group on $n$ letters and consider the lower central series $\Gamma^1_n = \mathcal{F}_n$, $\Gamma^{i+1}_n = [\Gamma^i_n,\mathcal{F}_n]$.
Then for each $q \geq 1$ there is a simplicial subset $B(q,G) \subseteq BG$, appearing in \cite{ACTG:CESSFCS}, whose set of $n$-simplices is $B_n(q,G) = \hom(\mathcal{F}_n / \Gamma^q_n, G) \subseteq G^{\times n}$. 
Thus an $n$-tuple in $G$ determines an $n$-simplex in $B(q,G)$ if and only if it generates a subgroup of $G$ of nilpotence class at most $q$.
The simplicial set $B(q,G)$ is always a partial group.
This is because the $L$ operation $BG_n = \hom(\mathcal{F}_n, G) \to \hom(\mathcal{F}_{2n},G) = BG_{2n}$ is given by the homomorphism $\mathcal{F}_{2n} \to \mathcal{F}_n$ sending $r_i$ to $r_{n-i+1}^{-1}$ for $i = 1, \dots, n$ and to $r_{i-n}$ for $i=n+1, \dots, 2n$.
As a homomorphism, this sends the subgroup $\Gamma^q_{2n}$ into $\Gamma^q_n$, so $B(q,G)$ is closed under the operation $L$ in $BG$, and hence is a partial group.
\end{example}

It will be helpful to have slightly more refined information about some partial compositions of the elements $L[f_1 | \cdots | f_n]$.
For the statement of the following lemma, we use the following notation for certain edges associated to $[f_1 | \cdots | f_n]$.
\begin{itemize}[left=0pt]
\item For $0 \leq i \leq n$, we write $f_{ii} = \id_{x_i} = s_0(x_i)$ where $x_i$ is the target of $f_i$ or the source of $f_{i+1}$. (Since $d_0f_i = d_1 f_{i+1}$ for $1\leq i \leq n-1$ there is no ambiguity here.)
\item For $i < j$ we write $f_{ij}$ for the total composite of $[f_{i+1} | \cdots | f_j]$.
\item For $j > i$ we write $f_{ji} \coloneqq f_{ij}^\dagger$, which is also the total composite of $[f_j^\dagger | \cdots | f_{i+1}^\dagger]$.
\end{itemize}
Notice that $f_{i-1,i} = f_i$.

\begin{lemma}\label{lem partial composites of L}
Suppose $X$ is a partial groupoid, $n > 0$, and $[f_1 | \cdots | f_n]$ is an $n$-simplex.
Let $[g_1 | \cdots | g_{2n}] = L[f_1 | \cdots | f_n]$.
Then for $i \leq j$ we have $g_{ij} = f_{|n-i|,|n-j|}$.
\end{lemma}
\begin{proof}
We have $g_t = f_{n+1-t}^\dagger \colon x_{n-t+1} \to x_{n-t}$ for $1 \leq t \leq n$ and $g_t = f_{t-n} \colon x_{t-n-1} \to x_{t-n}$ for $ n+1 \leq t \leq 2n$.
It follows that $f_{n-i,n-i} = \id_{x_{n-i}} = g_{ii}$ for $0\leq i \leq n$ and $f_{i-n,i-n} = \id_{x_{i-n}} = g_{ii}$ for $n \leq i \leq 2n$.
Now let's consider $g_{ij}$ for $i<j$.
This will be the total composite of one of the following:
\begin{align*}
[g_{i+1} | \cdots | g_j] &= [f_{n-i}^\dagger | \cdots | f_{n+1-j}^\dagger] && \text{for $i,j \leq n$} \\
[g_{i+1} | \cdots | g_j] &= [f_{i+1-n} | \cdots | f_{j-n}] && \text{for $i,j \geq n$} \\
[g_{i+1} | \cdots | g_j] &= [f_{n-i}^\dagger | \cdots | f_{j-n}] && \text{for $i <  n <  j$}
\end{align*}
Notice in the first two cases that we have $g_{ij} = f_{n-i,n-j}$ and $g_{ij} = f_{i-n,j-n}$, respectively.  

For the last case, set $k = n-i$ and $\ell = j-n$, so $[g_{i+1} | \cdots | g_j] = [f_k^\dagger | \cdots | f_\ell]$.
Generically write $\mu = \rho_{0,p}^*\colon X_p \to X_1$ for the total composite function. 
Suppose $k < \ell$.
Then
\begin{align*}
  g_{ij} = \mu [f_k^\dagger | \cdots | f_\ell] &= \mu([\mu ([f_k^\dagger | \cdots | f_k]) | f_{k+1} | \cdots | f_\ell ]) \\
  &= \mu([\id | f_{k+1} | \cdots | f_\ell ]) \\
  &= \mu([f_{k+1} | \cdots | f_\ell]) = f_{k\ell} = f_{n-i,j-n}.
\end{align*}
A similar calculation shows that $g_{ij} = f_{k\ell} = f_{n-i,j-n}$ when $k > \ell$.
We thus have $g_{ij} = f_{|n-i|,|n-j|}$ whenever $i\leq j$.
\end{proof}

\begin{lemma}\label{lem inv compat automatic}
Let $X$ and $Y$ be partial groupoids with associated anti-involutions $\nu$ and $\nu'$, and let $F\colon X \to Y$ a map of simplicial sets.
Then the square
\[ \begin{tikzcd}
X^\op \dar{\nu} \rar{F^\op} & Y^\op \dar{\nu'} \\
X \rar{F} &  Y 
\end{tikzcd} \]
commutes.
In particular, there is at most one anti-involution $\nu \colon X^\op \to X$ making a simplicial set $X$ into a partial groupoid.
\end{lemma}
\begin{proof}
The second statement follows by applying the first statement to  $F = \id_X$.
We now prove the first statement.
Let $F \colon X \to Y$ be a map of simplicial sets.
Since $Y$ is edgy, it suffices to show that the square commutes in simplicial degrees $0$ and $1$.
The first of these holds since $\nu$ and $\nu'$ are identities on vertices.
We write $(-)^\ddagger$ for $\nu_1' \colon Y_1 \to Y_1$.

Suppose $f\colon a\to b$ is in $X_1$. 
We consider the element $m\in Y_3$ obtained from $f^\dagger \in X_1$ via the following string.
\[ \begin{tikzcd}[row sep=tiny, column sep=small]
X_1 \rar["L"] & X_2 \rar{F} & Y_2 \rar{L'} & Y_4 \rar{d_0} & Y_3
 \\
{[f^\dagger]} \rar[mapsto] & {[f|f^\dagger]} \rar[mapsto] & {[Ff|Ff^\dagger]} \rar[mapsto] & {[F(f^\dagger)^\ddagger | (Ff)^\ddagger | Ff | F(f^\dagger)]} \rar[mapsto] & {[(Ff)^\ddagger | Ff | F(f^\dagger)]}
\end{tikzcd} \]
Since $d_1 [Ff | F(f^\dagger)] = Fd_1 [ f | f^\dagger] = F(\id_a) = \id_{Fa}$, we have $d_2(m) = [(Ff)^\ddagger | \id_{Fa}]$.
On the other hand, $d_1 [(Ff)^\ddagger | Ff ] = \id_{Fb}$, so we have $d_1(m) = [\id_{Fb} | F(f^\dagger)]$.
Both of these elements are degenerate.
We then have
\[
[(Ff)^\ddagger] = d_1 [(Ff)^\ddagger | \id_{Fa}] = 	d_1 d_2 (m) = d_1 d_1 (m) = d_1 [\id_{Fb} | F(f^\dagger)] = [F(f^\dagger)].
\]
This establishes that $F\nu_1 = \nu_1' F$, as required.
\end{proof}

In light of \cref{lem inv compat automatic}, we do not need to carry along the anti-involution as additional data, nor do we need to specify that maps are compatible with anti-involutions.
Maps will also be automatically compatible with the $L$ operators. 
Consequently, the following definition is reasonable. 
\begin{definition}\label{def pgpd category}
The category of partial groupoids, $\pgpd$, is the full subcategory of $\sset$ whose objects are partial groupoids.
\end{definition}

\section{Elements as matrices}
Suppose $X$ is an edgy simplicial set equipped with an anti-involution $\nu \colon X^\op \to X$. 
As usual, we write $(-)^\dagger \colon X_1 \to X_1$ for the action of $\nu$ on $X_1$.
We further assume that the map $X_0 \to X_0$ induced by $\nu$ is the identity, so in particular $\id_x^\dagger = \nu(s_0(x)) = s_0(\nu(x)) = s_0(x) = \id_x$ for each $x \in X_0$. 

\begin{definition}
Fix $n \geq 0$ and $\mathbf{f} \in X_n$. 
For $0\leq i \leq j \leq n$ we define
\begin{itemize}
\item $f_{ij}$ to be the edge $\rho_{ij}^*(\mathbf{f}) \in X_1$, and 
\item $f_{ji} \coloneqq f_{ij}^\dagger$.
\end{itemize}
The \emph{matrix form} of $\mathbf{f}$ is the $(n+1)\times (n+1)$ skew-symmetric matrix
\[
(f_{ij}) = 
\begin{bmatrix}
f_{00} & f_{01} & f_{02} & \cdots & f_{0n} \\
f_{10} & f_{11} & f_{12} & \cdots & f_{1n} \\
\vdots & & & & \vdots \\
f_{n0} & f_{n1} & f_{n2} & \cdots & f_{nn}
\end{bmatrix}
\]
in the involutive set $(X_1, \dagger)$.
\end{definition}
Isolating the superdiagonal gives $[f_{01} | f_{12} | \cdots | f_{n-1,n}]$, which is the usual way \eqref{standard form} we would write the element $\mathbf{f}$ when $n$ is positive.
Since $X$ is edgy, the function $X_n \to \mat_{n+1,n+1}(X_1)$ is injective (the case $n=0$ is just using that $s_0 \colon X_0 \to X_1$ is injective in any simplicial set, since $d_0$ is a left inverse).

\begin{lemma}\label{lem matrix form involution}
Using the matrix form, $\nu_n \colon X_n \to X_n$ acts by transposition.
\end{lemma}
\begin{proof}
Suppose $i\leq j$, and $\mathbf{g} \in (X^\op)_n$.
Then the $j,i$ entry in the matrix associated to $\nu(\mathbf{g}) = \mathbf{f} \in X_n$ is given by
\[
	f_{ji} = f_{ij}^\dagger = [\rho_{ij}^*(\nu(\mathbf{g}))]^\dagger = [\nu(\rho_{ij}^*(\mathbf{g}))]^\dagger = \rho_{ij}^*(\mathbf{g}) = g_{ij}.
\]
Since both matrices are skew-symmetric, this implies that $(f_{ij})^T = (g_{ij})$.
\end{proof}

\begin{lemma}\label{lem matrix form order preserving}
If $\alpha \colon [m] \to [n]$ is an order-preserving map and $(f_{ij})$ is the matrix form of $\mathbf{f} \in X_n$, then $(f_{\alpha k, \alpha \ell})$ is skew-symmetric and the matrix form for $\alpha^* \mathbf{f} \in X_m$.
\end{lemma}
\begin{proof}
For $k \leq \ell$ in $[m]$, we have 
$\rho_{k,\ell}^* \alpha^* \mathbf{f} = (\alpha \rho_{k,\ell})^* \mathbf{f} = \rho_{\alpha k, \alpha \ell}^* \mathbf{f} = f_{\alpha k, \alpha \ell}$.
Of course $(f_{\alpha k, \alpha \ell})$ is a skew-symmetric matrix since $(f_{ij})$ is, so it must be the matrix form for $\alpha^*\mathbf{f}$.
\end{proof}

\begin{definition}\label{def germinable}
Suppose $X$ is an edgy simplicial set with anti-involution, and $\mathbf{f} \in X_n$. 
Write $\chi = \chi_n \colon [2n] \to [n]$ for the function with $\chi(i) = |n-i|$.
We say that $\mathbf{f}$ is \emph{germinable} if there is $L\mathbf{f} \in X_{2n}$ whose matrix form is given by $(f_{\chi i, \chi j})$.
\end{definition}

Notice that the matrix $(f_{\chi i, \chi j})$ is automatically skew-symmetric since $(f_{ij})$ is (though it is not necessarily the matrix form for any element of $X$).
As an example, when $n=2$ the matrix for $L\mathbf{f}$ is given by
\[
\begin{bmatrix}
f_{00} & f_{01} & f_{02} \\
f_{10} & f_{11} & f_{12} \\
f_{20} & f_{21} & f_{22}
\end{bmatrix}
\mapsto
\begin{bmatrix}
f_{22} & f_{21} & f_{20} & f_{21} & f_{22} \\
f_{12} & f_{11} & f_{10} & f_{11} & f_{12} \\
f_{02} & f_{01} & f_{00} & f_{01} & f_{02} \\
f_{12} & f_{11} & f_{10} & f_{11} & f_{12} \\
f_{22} & f_{21} & f_{20} & f_{21} & f_{22}
\end{bmatrix}
\]
This matrix consists of the matrices for $\nu\mathbf{f}$ and $\mathbf{f}$ glued together along the central $f_{00}$.
The diagonal and anti-diagonal consist of identities.
The $3\times 3$ matrices in the northeast and southwest are also interesting and related to (reflected about central column and central row) the original matrix. This behavior is typical for all $n$.
The matrices for $L\mathbf{f}$ are always centrosymmetric (see \cite{Saibel:NICM,Good:ICM}) since
\[
	g_{2n-k,2n-\ell} = f_{|n-(2n-k)|,|n-(2n-\ell)|} = f_{|k-n|, |\ell-n|} = f_{|n-k|, |n-\ell|} = g_{k,\ell}
\]
by \cref{lem partial composites of L}.

\begin{proposition}\label{prop germinable}
Let $X$ be an edgy simplicial set with anti-involution.
Then every element of $X$ is germinable if and only if $X$ is a partial groupoid.
\end{proposition}
\begin{proof}
Let $\mathbf{f}$ be represented by $[f_{01} | \cdots | f_{n-1,n}] = [f_1 | \cdots | f_n]$, and suppose that $\mathbf{f}$ is germinable. 
Then the superdiagonal of $(f_{\chi i, \chi j})$ is \[ [f_{n-1,n}^\dagger | \cdots | f_{01}^\dagger | f_{01} | \cdots | f_{n-1,n}] =  [f_n^\dagger| \cdots | f_1^\dagger | f_1 | \cdots | f_n],\] which represents an element of $X_{2n}$.
The total composite $\rho_{0,2n}^*L\mathbf{f} = f_{\chi(0), \chi(2n)} = f_{nn}$ is an identity.
Thus if every element of $X$ is germinable, then $X$ is a partial groupoid.

Conversely, suppose $\mathbf{f}$ is an element of a partial groupoid and $\mathbf{g} = L\mathbf{f}$. 
By \cref{lem partial composites of L}, we have $g_{ij} = f_{\chi i, \chi j}$ for $i \leq j$.
Since $(g_{ij})$ and $(f_{\chi i, \chi j})$ are both skew-symmetric, they must be equal. 
Hence $\mathbf{f}$ is germinable.
\end{proof}

Recall that $J\colon \simpcat \to \fplus$ denotes the subcategory inclusion, which induces a faithful restriction functor $J^* \colon \symset \to \sset$.

\begin{proposition}\label{prop spiny to pgpd}
Suppose $Z \colon \fplus^\op \to \set$ is a symmetric set.
If $J^*Z$ is edgy, then it is a partial groupoid.
\end{proposition}
\begin{proof}
There is a canonical anti-involution on $X = J^*Z$ given by $\nu_n(x) = \tau_n^*x$ for $x\in X_n = Z_n$. 
If $\alpha \colon [n] \to [m]$ is an order-preserving map and $x\in X_m$, then \[ \alpha^* \nu_m(x) = \alpha^* \tau_m^*(x) = (\tau_m \alpha)^*(x) = (\tau_m \alpha \tau_n \tau_n)^*(x) = \tau_n^* (\tau_m \alpha \tau_n)^* x = \nu_n(\alpha^*x)\]
where the $\alpha^*$ on the far right is occurring in $(J^*Z)^\op$ and the one on the far left is occurring in $J^*Z$.
This shows that $\nu$ constitutes a map of simplicial sets $\nu \colon X^\op \to X$.
Since $\tau_n$ is an involution, so too is $\nu_n$; when $n=0$ they are both identities.

In light of \cref{prop germinable}, we need to show that every element of $X$ is germinable. 
If $\mathbf{f} \in X_n = Z_n$, then $f_{ij} = \rho_{ij}^*(\mathbf{f})$ for any $i,j \in [n]$, including when $i > j$ (in which case $\rho_{ij}$ is not order-preserving).
To see this, note that for $i\leq j$ we have \[ f_{ji} = f_{ij}^\dagger = (\rho_{ij}^*(\mathbf{f}))^\dagger
= 
\tau_1^* \rho_{ij}^* \mathbf{f} = \rho_{ji}^*\mathbf{f}
\]
using the definition of $\nu_1$.
We let $L\mathbf{f} \coloneqq \chi_n^*(\mathbf{f}) \in Z_{2n}= X_{2n}$.
The matrix form for $L\mathbf{f}$ is precisely $(f_{\chi i, \chi j})$, hence $\mathbf{f}$ is germinable.
\end{proof}

\section{Spiny symmetric sets}\label{sec spiny}
We just saw that if the underlying simplicial set of a symmetric set is edgy, then the underlying simplicial set is a partial groupoid.
The present section gives an alternative characterization of the symmetric sets with edgy restriction, as the `spiny symmetric sets.'
Later sections do not logically depend upon this one, and the reader may take the statement of \cref{prop spiny edgy} as the definition of `spiny' and move on to the next section.

\begin{definition}\label{def spines}
Let $S$ be a finite set with at least two elements. 
A \emph{spine} for $S$ is a tree $T$ with $S$ as its vertex set, i.e.\ a spanning tree for the complete graph on $S$.
We will also present a spine as an ordered list of two element subsets $e_1, \dots, e_{|S|-1}$ of the set $S$ (the edges of the tree).
A \emph{system of spines} consists of a chosen spine on each of the sets $[n]$ for $n\geq 1$.
\end{definition}

Of course there is only a single spine for a 2-element set, and the three spines on a 3-element set are isomorphic. 
Cayley's theorem \cite{Cayley:TOT} says there are $|S|^{|S|-2}$ spines on $S$.

\begin{example}\label{ex systems}
Nice examples of systems of spines include the following:
\begin{enumerate}
\item $e_i = \{i-1,i\}$ for $1\leq i \leq n$ --- we call this the `standard system of spines' \label{usual spine system}
\item $e_i = \{0,i\}$ for $1\leq i \leq n$
\item $e_i = \{i-1,n\}$ for $1\leq i \leq n$
\end{enumerate}
The latter two are essentially the same spine, while the first is different from them when $n\geq 3$.
\end{example}

The standard system of spines \eqref{usual spine system} is at the heart of the Grothendieck--Segal condition, while the second system of spines underlies the Bousfield--Segal condition \cite{Bergner:AID2,HackneyBergner:GASO,Stenzel:BSS}.

Recall that $\finplus$ denotes the category of finite nonempty sets, and the inclusion $\fplus \to \finplus$ is an equivalence of categories, hence the restriction $[\finplus^\op, \set] \to [\fplus^\op, \set] = \symset$ is also an equivalence.
We temporarily work with this more general notion of symmetric set.

\begin{definition} Let $X\colon \finplus^\op \to \set$ be a presheaf.
\begin{itemize}
\item Given a spine $T = e_1, \dots, e_{|S|-1}$ for a set $S$, we say that $X$ is \emph{injective along $T$} if 
\[
	\mathscr{E}_T \colon X(S) \to \prod_{i=1}^{|S|-1} X(e_i)
\]
is injective. (Here, the maps $e_i \to S$ are the subset inclusions.)
\item Given a system of spines $E = \{ T_1, T_2, \dots, T_n, \dots \}$, we say that $X$ is \emph{injective along $E$} if $X$ is injective along $T_n$ for each chosen spine $T_n$ of $[n]$ in the collection $E$.
\end{itemize}
\end{definition}

The second bullet in this definition works equally well for presheaves on $\fplus$, using the order preserving maps $[1] \to [n]$ with image $e_i$, yielding maps of the form
\[
	X_n \to \prod_{i=1}^n X_1
\]
which we ask to be injections (and also call $\mathscr{E}_T$).

\begin{example}
Let $E$ be the standard system of spines from \cref{ex systems}\eqref{usual spine system}.
Then a symmetric set $X$ is injective along $E$ if and only if $J^*X$ is an edgy simplicial set, since if $T \in E$ is the spine on $[n]$, then $\edgemap_T$ is the map  $\edgemap_n$ appearing in \cref{def edgy sset}.
\end{example}

\begin{lemma}\label{spine permuting lemma}
Suppose $T$ is a spine of $S$ and $X$ is injective along $T$.
Let $\sigma \colon S \to S'$ be a bijection, and declare $T'$ to be the spine of $S'$ obtained from $T$ along this bijection.
Then $X$ is injective along $T'$.
\end{lemma}
\begin{proof}
Let $\sigma_i \colon e_i \to e_i' \coloneqq \sigma(e_i)$ be the restriction of $\sigma$ to the two element subset $e_i \subseteq S$.
The $e_i'$ form the edges of $T'$.
The following diagram commutes,
\[ \begin{tikzcd}
X(S') \dar["\sigma^*"', "\cong"] \rar{\mathscr{E}_{T'}} & \prod X(e_i')  \dar["\cong"', "\prod(\sigma_i)^*"] \\
X(S) \rar{\mathscr{E}_T} & \prod X(e_i)
\end{tikzcd} \]
and the bottom map was assumed injective. Hence the top map is injective.
\end{proof}

\begin{theorem}\label{theorem spines system}
Suppose $X$ is injective along $E$ for some system of spines $E$.
Then $X$ is injective along $T$ for every spine $T$.
\end{theorem}
\begin{proof}
By \cref{spine permuting lemma}, it is enough to prove the result for spines on $[n]$, for all $n\geq 1$.
We prove the result by induction on $n$.
When $n=1$, there is a single spine $T$ and $\mathscr{E}_T$ is the identity on $X_1$, which is an injective map.
Let $n > 1$, and suppose that $X$ is injective along $T'$ for every spine $T'$ of $[n-1]$.

Let $T$ be an arbitrary spine of $[n]$, and let $T_n$ be the spine of $[n]$ that is in the system $E$.
As $T$ is a tree on a set with size at least two, there is a leaf vertex $v$, which is connected to a single other vertex $w$.
By permuting the elements of $[n]$, we can suppose that $v = n$ and $w = n-1$, and by \cref{spine permuting lemma} it suffices to check injectivity on this alternate spine.
The same can be done for $T_n$.
We thus assume that both $T$ and $T_n$ have $n$ has a leaf vertex attached only to $n-1$.
Consider the following commutative diagram.
\[ \begin{tikzcd}
& X_n \ar[dl,"\mathscr{E}_T"'] \ar[dr,"\mathscr{E}_{T_n}", hook] \dar["d_n \times d_0^{n-1}" description] \\
\left( \prod\limits_{i=1}^{n-1} X_1\right) \times X_1 & X_{n-1} \times X_1  \lar[hook'] \rar & \prod\limits_{i=1}^n X_1 
\end{tikzcd} \]
By assumption (perhaps invoking \cref{spine permuting lemma}) we know that $\mathscr{E}_{T_n}$ is injective, which implies that the vertical map $X_n \to X_{n-1} \times X_1$ is also injective.
By the inductive hypothesis, \[ \mathscr{E}_{T|_{[n-1]}} \colon X_{n-1} \to \prod_{i=1}^{n-1} X_1\]
is injective, so we conclude that $\mathscr{E}_T$ is injective, as desired.
\end{proof}

\begin{definition}\label{def spiny symset}
Let $X$ be a symmetric set.
Then $X$ is said to be \emph{spiny} if it is injective along $E$ for some system of spines $E$.
\end{definition}
In light of \cref{theorem spines system}, the word `some' may be replaced by the word `every' in the preceding definition.

\begin{proposition}\label{prop spiny edgy}
A symmetric set $X$ is spiny if and only if $J^*X$ is an edgy simplicial set.
\end{proposition}

\section{Partial groupoids as spiny symmetric sets}

Our next goal is to show that the category of partial groupoids (\cref{def pgpd category}) is equivalent to the category of spiny symmetric sets. 
In this section, symmetric set will always mean a presheaf for the skeletal category $\fplus$, and this equivalence will turn out to be an isomorphism of categories (\cref{thm spiny equals pgpd}).
This isomorphism is obtained by restricting along the functor $J\colon \simpcat \to \fplus$.
Eventually we will prove the following, which gives a bijection at the level of objects between the two categories.

\begin{theorem}\label{thm symm nerve}
If $X$ is a partial groupoid, then there is a unique symmetric set $Z$ with $J^*Z = X$.
\end{theorem}
The equality implies, in particular, that $Z_n = X_n$ for all $n$.
By \cref{prop spiny edgy} we know that $Z$ is spiny in this case.

To prove existence of $Z$, we wish to extend the action of $\simpcat$ on $X$ to an action of $\fplus$, which we do using the (skew-symmetric) matrix forms for elements.
Namely, we want to define $\phi^*(f_{ij})$ to be $(f_{\phi i, \phi j})$ whenever $(f_{ij})$ is the matrix form for an element of $X_n$ and $\phi \colon [m] \to [n]$ is a function.
According to \cref{lem matrix form order preserving}, this would be an extension of the existing $\simpcat$-action.
But it must be shown that this formula is well-defined for maps which are not order-preserving, and we do this in \cref{prop Z existence}.

We currently have precisely two classes of non-order preserving maps $\phi$ where we can guarantee that $(f_{\phi i, \phi j})$ is the matrix form for some element of $X$. 
The first class consists of the flip maps $\tau_n$.
In this case, the matrix for $\tau_n^*\mathbf{f}$ is the transposition of $(f_{ij})$, which is the matrix form for $\nu_n\mathbf{f}$ by \cref{lem matrix form involution}.

The second class consists of the fold functions $\chi_n \colon [2n] \to [n]$ with $\chi_n(i) = |n-i|$.
If $\mathbf{f} \in X_n$, then since $\mathbf{f}$ is germinable, the element $L\mathbf{f} \in X_{2n}$ has the correct matrix form for $\chi_n^*\mathbf{f}$ by \cref{def germinable}.
Fortunately, the morphisms of $\fplus$ are generated by the maps $\chi_n$ and the order-preserving maps.

\begin{lemma}\label{lem iterated fold} % See Jan 26 note
The composite 
\[ \chi_n^w \coloneqq \chi_n \circ \chi_{2n} \circ \chi_{2^2n}\circ  \dots \circ \chi_{2^{w-1} n} \colon [2^w n] \to [2^{w-1}n] \to \cdots \to [2n] \to [n] \] is given, for $0\leq t \leq 2^{w-1} - 1$ and $2tn \leq i \leq 2(t+1)n$ by
\begin{equation}\label{eq chi n w formula}
\chi_n^w(i) = |(2t+1)n-i| = 	\begin{cases}
		(2t+1)n-i & i \in [2tn, (2t+1)n] \\
		i-(2t+1)n & i \in [(2t+1)n, (2t+2)n].
	\end{cases}
\end{equation}
\end{lemma}
\begin{proof}
For the remainder of the proof we take $\chi_n^w$ to mean the formula given in \eqref{eq chi n w formula}; it is immediate that $\chi_n^1 = \chi_n \colon [2n] \to [n]$.
We will calculate $\chi_n^w \circ \chi_{2^wn} \colon [2^{w+1}n] \to [n]$ and show that it is equal to $\chi_n^{w+1}$, which will complete the proof.
For a fixed $i \in [2^{w+1}n]$, we have $\chi_n^w \chi_{2^wn}(i) = \chi_n^w(|2^wn - i|),$ and we split into cases based on the absolute value.
Let $0\leq t \leq 2^w-1$ be such that $2tn \leq i \leq (2t+2)n$.
\begin{itemize}[left = 0pt]
\item 
If $t \leq 2^{w-1}-1$, then $i\leq 2^wn$, so we have $\chi_{2^wn}(i) = 2^wn - i$.
Setting $s = 2^{w-1} - (t+1)$, we have
\[
	2sn = 2^wn - (2t+2)n \leq 2^wn - i \leq 2^wn - 2tn = 2(s+1)n.
\]
By \eqref{eq chi n w formula} we have
\begin{align*}
	\chi_n^w(2^wn - i) = |(2s+1)n-(2^wn-i)|  &= | (2^w - 2t - 1)n - (2^wn - i)| \\
	&= |i - (2t+1)n| = \chi_n^{w+1}(i).
\end{align*}
\item If $t \geq 2^{w-1}$, then $i \geq 2^wn$, so we have $\chi_{2^wn}(i) = i - 2^wn$.
Setting $s=t-2^{w-1}$, we have
\[
	2sn = 2tn - 2^wn \leq i - 2^wn \leq 2(t+1)n - 2^wn = 2(s+1)n.
\]
By \eqref{eq chi n w formula} we have
\begin{align*}
	\chi_n^w(i-2^wn) = |(2s+1)n-(i-2^wn)|  &= 
	 | (2t-2^w +1)n - (i-2^wn)| \\
	&= |(2t+1)n-i| = \chi_n^{w+1}(i).
\end{align*}
\end{itemize}
Thus in both cases we have $\chi_n^w \chi_{2^wn}(i) = \chi_n^{w+1}(i)$.
\end{proof}

\begin{proposition}\label{prop fplus generation}
The category $\fplus$ is generated by the morphisms in $\simpcat$ along with the fold maps $\chi_n \colon [2n] \to [n]$.
If $\phi \colon [m] \to [n]$ is an arbitrary function, then we can factor $\phi$ as
\[
	[m] \xrightarrow{\alpha} [2^wn] \xrightarrow{\chi_n^w} [n]
\]
where $\alpha$ is order-preserving. 
\end{proposition}
\begin{proof}
Let $w \geq \log_2(m+1)+1$.
Define $\alpha(k) \coloneqq \phi(k) + (2k+1)n$, which is order-preserving: since $\phi$ lands in $[n]$, we have \[ (2k+1)n \leq \alpha(k) \leq (2k+2)n \leq (2(k+1)+1)n \leq \alpha(k+1) \leq (2(k+1)+2)n.\]
Notice that $\alpha(m) \leq 2(m+1)n \leq 2 \cdot 2^{w-1} n$, so $\alpha$ indeed has the indicated codomain.
Using the bounds above and the formula from \cref{lem iterated fold}, 
\[
	\chi_n^w(\alpha (k)) = \chi_n^w (\phi(k) + (2k+1)n) = \phi(k) + (2k+1)n - (2k+1)n = \phi(k). \qedhere
\]
\end{proof}

\begin{proposition}\label{prop Z existence}
If $X$ is a partial groupoid, then for each $\mathbf{f} \in X_n$ having matrix form $(f_{ij})$ and each morphism $\phi \colon [m] \to [n]$ of $\fplus$, there is a unique element $\phi^*\mathbf{f} \in X_m$ having matrix form $(f_{\phi i, \phi j})$.
This defines a symmetric set $Z$ extending $X$.
\end{proposition}
\begin{proof}
Suppose $\mathbf{f} \in X_n$, and $\phi \colon [m] \to [n]$ is an arbitrary function.
Since $X_m \to \mat_{m+1,m+1}(X_1)$ is injective, there is at most one such element $\phi^*\mathbf{f}$.
Factor $\phi$ as $\chi \circ \alpha$ as in \cref{prop fplus generation}, where $\alpha$ is order-preserving and $\chi$ is a composite of fold maps $\chi_k$.
Since every element of $X$ is germinable in the sense of \cref{def germinable}, $(f_{\chi i, \chi j})$ is the matrix for some element of $X$.
By \cref{lem matrix form order preserving}, we also know that $(f_{\chi \alpha i, \chi \alpha j})$ is the matrix for an element of $X$.
But this matrix is just $(f_{\phi i, \phi j})$, so $\phi^*\mathbf{f} \in X_m$ is well-defined.

To show that we have a functor $\fplus^\op \to \set$, we must show that $\psi^*\phi^* \mathbf{f} = (\phi \psi)^* \mathbf{f}$ and $(\id)^* \mathbf{f} = \mathbf{f}$.
The latter is immediate.
For a composition $\phi \circ \psi$ we have
\[
  \psi^* \phi^* (f_{ij}) = \psi^*(f_{\phi i, \phi j}) = (f_{\phi \psi i, \phi \psi j}) = (\phi \psi)^*(f_{ij}).
\]
We thus have a functor $Z \colon \fplus^\op \to \set$ with $Z_n = X_n$.
By \cref{lem matrix form order preserving}, we know that $J^*Z = X$, that is, this $\fplus$-action on $Z$ extends the usual $\simpcat$-action on $X$.
\end{proof}

\begin{lemma}\label{lem fully faithful}
The functor $J^* \colon \symset \to \sset$ restricts to a fully-faithful functor from spiny symmetric sets to partial groupoids.
\end{lemma}
\begin{proof}
The fact that the functor restricts as indicated is \cref{prop spiny to pgpd} and \cref{prop spiny edgy}.
This restriction is automatically faithful since $J \colon \simpcat \to \fplus$ is the identity on objects.
Let $Y,Z$ be spiny symmetric sets and suppose $F \colon J^*Y \to J^*Z$ is a map of simplicial sets.
By the proof of \cref{prop spiny to pgpd}, the anti-involutions on these partial groupoids are given by $\tau_n^*$ and the maps $L$ from \cref{def partial groupoids} are given by $\chi_n^*$.
Thus, by \cref{lem inv compat automatic}, we have that $\tau_n^* F_n = F_n \tau_n^*$.
In particular this holds for $n=1$, so we see for $[f_1 | \cdots | f_n] \in Y_n$ that
\begin{align*}
	F\chi_n^*[f_1 | \cdots | f_n] = F[f_n^\dagger | \cdots | f_1^\dagger | f_1 | \cdots | f_n] &= [F(f_n^\dagger) | \cdots | F(f_1^\dagger) | Ff_1 | \cdots | Ff_n] \\
	&=[(Ff_n)^\ddagger) | \cdots | (Ff_1)^\ddagger) | Ff_1 | \cdots | Ff_n]
\end{align*}
is equal to $\chi_n^* F[f_1 | \cdots |f_n] \in Z_{2n}$.
Thus $\chi_n^* F_n = F_{2n} \chi_n^*$.
By \cref{prop fplus generation} the morphisms in $\fplus$ are generated by the order-preserving maps and the fold maps $\chi_n$, hence $\phi^* F_n = F_m \phi^*$ for any function $\phi \colon [m] \to [n]$.
We conclude that $J^*$ restricted to the subcategory of spiny symmetric sets is full.
\end{proof}

With more work, one can show that you only need the codomain to be a spiny symmetric set in order to lift maps.
We will not need this strengthened result.

\begin{proof}[Proof of \cref{thm symm nerve}]
It is shown in \cref{prop Z existence} that there exists a symmetric set $Z$ with $J^*Z = X$.
Suppose $Z' \colon \fplus^\op \to \set$ is a functor with $J^*Z' = X$.
By \cref{lem fully faithful}, the identity $J^*Z \to J^*Z'$ lifts to a map of symmetric sets $Z\to Z'$.
Since this is the identity at each level, it is an identity of symmetric sets.
\end{proof}

\begin{theorem}\label{thm spiny equals pgpd}
The category of spiny symmetric sets is isomorphic to the category of partial groupoids.
\end{theorem}
\begin{proof}
The functor $J^*$ from spiny symmetric sets to partial groupoids is bijective on objects by \cref{thm symm nerve} and fully-faithful by \cref{lem fully faithful}.
\end{proof}

Since a symmetric set $Z$ is reduced if and only if $J^*Z$ is reduced, \cref{gonzalez theorem} and \cref{thm spiny equals pgpd} yield the following.
\begin{corollary}\label{cor spiny partial groups}
The category of reduced spiny symmetric sets is equivalent to the category of partial groups. \qed
\end{corollary}

\section{Reflections onto reduced and spiny symmetric sets}
\label{section reflections}
We have a square of fully faithful functors.
\[ \begin{tikzcd}
\pgrp \rar \dar & \symset_* \dar \\
\pgpd \rar & \symset
\end{tikzcd} \]
For convenience, we regard these as honest subcategories for the remainder of the paper, so partial group is synonymous with reduced spiny symmetric set, and similarly for partial groupoids.
Our goal in this section is to show that partial groupoids and partial groups form \emph{reflective} subcategories of the category of symmetric sets.

A full subcategory $D \subseteq C$ is a \emph{reflective subcategory} just when the inclusion functor has a left adjoint $C \to D$ (see, for instance, \cite[Definition 4.5.12]{Riehl:CTC}).
To specify this left adjoint is to give, for each object $c\in C$, a \emph{reflection arrow} $\eta_c \colon c\to \bar{c}$ where $\bar{c} \in D$, satisfying the universal property that any map $c\to d$ with codomain in $D$ factors as $f\circ \eta_c$ for a unique $f\colon \bar{c} \to d$ in $D$. 
In other words, the map $\hom(\bar{c}, d) \to \hom(c,d)$ given by precomposition with $\eta_c$ is a bijection for each $d\in D$.
The left adjoint is called the \emph{reflection} or \emph{reflector}.

There are various criteria in the literature for establishing when a subcategory is reflective.
A number of these will suffice to establish that the subcategories in the above square are reflective.
The proofs in this section follow the direct path, by constructing explicit left adjoints to the inclusion functors.
This choice was made to ensure the proofs are accessible and rely only on elementary category theory.
We also anticipate that the explicit reflections will be useful when computing colimits of partial groups and partial groupoids (\S\ref{sec consequences}).
However, we briefly indicate in \cref{rmk formal spinification} and \cref{remark alt proof} two alternative ways to arrive at this fact, relying on general theorems in category theory.
We recommend \cref{remark alt proof} for those with the appropriate background, and it may be read now, in place of other proofs in this section.

\subsection{The reduction of a symmetric set}
Given $X \in \symset$, we write $X_0$ for the constant symmetric set with values $X_0([n]) = X_0$ and on which each morphism in $\fplus$ acts as the identity. 
The symbol $\ast$ denotes the terminal symmetric set. 
\begin{definition}\label{def reduction}
If $X$ is a symmetric set, let $\mathcal{R}X$ be its \emph{reduction}, defined via the following pushout diagram
\[ \begin{tikzcd}
X_0 \rar[hook] \dar \ar[dr, phantom, "\ulcorner" very near end] & X \dar \\
\ast \rar & \mathcal{R}X.
\end{tikzcd} \]
The monomorphism $X_0 \hookrightarrow X$ is induced from the unique maps $[n]\to[0]$. 
\end{definition}

More concretely, an element $x$ in $X_n$ is \emph{fully degenerate} if there exists an element $x'$ in $X_0$ such that $x = \psi^*x'$, where $\psi \colon [n] \to [0]$ is the unique map.
Notice that if $x\in X_n$ is fully degenerate and $\phi \colon [m] \to [n]$ is any function, then $\phi^*x \in X_m$ is also fully degenerate. 
So long as $X$ is non-empty, $\mathcal{R}X_n$ is the quotient of $X_n$ where we have identified the set of fully degenerate elements, and the usual action of $\fplus$ descends to this quotient. (If $X$ is empty, then $\mathcal{R}X = \ast$.)
This symmetric set is indeed reduced ($\mathcal{R}X_0$ is a point), and there is a canonical quotient map $X \to \mathcal{R}X$ given by identifying the fully degenerate elements at each level.
The following is immediate from \cref{def reduction}.

\begin{lemma}
If $Y$ is a reduced symmetric set (i.e.\ $Y_0$ is a point) then any map $F \colon X \to Y$ factors uniquely as $X \to \mathcal{R}X \to Y$. \qed
\end{lemma}

This proves that the category of reduced symmetric sets is a reflective subcategory of symmetric sets.

\begin{proposition}
If $X$ is a spiny symmetric set, then so is its reduction $\mathcal{R}X$.
\end{proposition}
\begin{proof}
If $x = [f_1|\cdots | f_n]$ and $y = [g_1 | \cdots | g_n]$ are elements of
$X_n$ whose images $\bar{x}, \bar{y} \in (\mathcal{R}X)_n$ have the same edge
profile, then for each $i$, either $f_i = g_i$ or both $f_i$ and $g_i$ are 
degenerate (i.e.\ $f_i$ and $g_i$ are identities, possibly on different objects).  
If $f_i$ and $g_i$ are degenerate for each $i$, then $x$ and $y$ are themselves 
fully degenerate, and hence $\bar{x} = \bar{y}$. Otherwise, fixing $k$ with
$f_k = g_k$, move outward to the boundary of the words from position $k$ to see
that $f_i$ and $g_i$ have the same source and target for each $i$, hence are
equal even when both are degenerate. So $x = y$ in this case. 
\end{proof}

\begin{example}[Word classifiers]\label{example reduced representable}
A chaotic (or indiscrete) groupoid is a category having a unique morphism between any two objects.
The representable presheaf $\rep{m} = \hom(-,[m])$ on $\fplus$ is the nerve of the chaotic groupoid on $m+1$ objects, hence is a spiny symmetric set.
Its reduction $\fpg{m} \coloneqq \mathcal{R}\rep{m}$ has the property that for any other reduced symmetric set $X$ we have (using adjointness and the Yoneda lemma)
\[
	\hom(\fpg{m}, X) = \hom(\mathcal{R}\rep{m}, X) \cong \hom(\rep{m}, X) \cong X_m.
\]
Thus the reduced spiny symmetric set $\fpg{m}$ can be considered the free partial group on a length $m$ admissible sequence.
Since $\rep{m}_n = \hom([n],[m])$ has $(m+1)^{n+1}$ elements, 
and exactly $m+1$ of these are fully degenerate, we see that $\fpg{m}_n$ has cardinality $(m+1)^{n+1} - m$.

Let's briefly consider the case $m=2$.
In degree 1, we have nondegenerate edges 
\begin{align*}
a_1 &\colon 0 \to 1 & a_1^\dagger &\colon 1 \to 0 \\
a_2 &\colon 1\to 2 & a_2^\dagger &\colon 2 \to 1 \\
a_2 \circ a_1 = a_3 &\colon 0 \to 2 & a_1^\dagger \circ a_2^\dagger = a_3^\dagger &\colon 2 \to 0
\end{align*}
along with the three degenerate edges in $\rep{2}$ or one degenerate edge in $\fpg{2} = \mathcal{R}\rep{2}$.
\end{example}

\begin{example}[Reduction of groupoids]\label{ex reduction groupoid}
Any groupoid $\mathscr{G}$ naturally determines a partial group.
The elements of the partial group consist of non-identity morphisms of $\mathscr{G}$, along with a single unit element $e$.
A word in this set is admissible just when the corresponding word with units $e$ thrown away is composable in $\mathscr{G}$.
An alternative description is to take the nerve of $\mathscr{G}$, which is a spiny symmetric set, and then apply $\mathcal{R}$.
\end{example}

\subsection{The reflection onto spiny symmetric sets}
Our next aim is to define a reflection from the category of symmetric sets to the category of spiny symmetric sets.
This is more delicate than the reduction of a symmetric set.
For a symmetric set $X$ and $n\geq 1$, recall the map $\edgemap_n \colon X_n \to \prod_{i=1}^n X_1$ from \cref{def edgy sset}, and let $\edgemap_0 \colon X_0 \to X_0$ be the identity.

\begin{definition}\label{def flat relation} 
For a symmetric set $X$ and $x, x' \in X_n$, we write $x \approx x'$ whenever there exists a map $\phi \colon [n] \to [m]$ and elements $\tilde x, \tilde x' \in X_m$ with
\begin{enumerate}
\item $\edgemap_m(\tilde x) = \edgemap_m(\tilde x')$, (if $m=0$ this means $\tilde x = \tilde x'$) \label{item E}
\item $\phi^*\tilde x = x$, and \label{item x}
\item $\phi^* \tilde x' = x'$. \label{item x prime}
\end{enumerate}
Then $\approx$ is a reflexive and symmetric relation on $X_n$. 
Let $\sim$ be the equivalence relation generated by $\approx$, and write $[x]$ for the equivalence class of $x$. 
\end{definition}

If $\psi \colon [k] \to [n]$ is any map and $\tilde x ,\tilde x'  \in X_m$ and $\phi \colon [n] \to [m]$ are witnesses for $x \approx x' \in X_n$ as above, then $\tilde x, \tilde x' $ and $\phi \psi$ are witnesses for $\psi^*x \approx \psi^*x'$ in $X_k$. 
In this way one sees that the action of $\fplus^\op$ on $X$ preserves the equivalence relation $\sim$. 

\begin{definition}\label{def flat simplicial set}
For a symmetric set $X$, we write $X^\flat$ for the quotient of $X$ by $\sim$. 
\end{definition}
Thus, if $\psi \colon [k] \to [n]$ is a function and $[x] \in X_n^\flat$, then $\psi^*[x] \coloneqq [\psi^*x] \in X^\flat_k$.
The object $X^\flat \in \symset$ may also be described as pushout of symmetric sets, using the maps $c_n$ from \cref{remark alt proof}.

Notice that if $X$ is a spiny symmetric set, the relation $\sim$ is the identity relation, so the natural map of symmetric sets $X \twoheadrightarrow X^\flat$ is an isomorphism. 
Also notice that $X_0 \to X^\flat_0$ is always a bijection.

\begin{lemma}\label{lem flat}
If $U$ is a spiny symmetric set, then $X \twoheadrightarrow X^\flat$ induces a bijection $\hom(X^\flat, U) \to \hom(X,U)$.
\end{lemma}
\begin{proof}
Since $X \twoheadrightarrow X^\flat$ is an epimorphism, the indicated map is a monomorphism without any assumption on $U$.
Let $F\colon X \to U$ be a map of symmetric sets.
Suppose $x \approx x'$ in $X_n$, witnessed by $\phi, \tilde x, \tilde x'$ as above.
If $m=0$ then since $\edgemap_0 \colon X_0 \to X_0$ is the identity, we have $\tilde x=\tilde x'$.
In any other case $m\geq 1$, we have the square
\[ \begin{tikzcd}
X_m \dar{F} \rar{\edgemap_m} & \prod X_1 \dar{F} \\
U_m \rar[hook, "\edgemap_m"] & \prod U_1
\end{tikzcd} \]
Since the top map identifies $\tilde x$ and $\tilde x'$ and the bottom map is an injection, we have $F(\tilde x) = F(\tilde x')$.
Since $F(x) = F(\phi^* \tilde x) = \phi^*F(\tilde x)$, and similarly for $x'$, we conclude that $F(x) = F(x')$.
Thus $F \colon X \to U$ factors uniquely as $X \to X^\flat \to U$, and we conclude that $\hom(X^\flat, U) \to \hom(X,U)$ is a bijection.
\end{proof}

One observation is that if $x,x' \in X_n$ are such that $\edgemap_n(x) = \edgemap_n(x')$, then $[x]=[x'] \in X_n^\flat$.
We have made some progress regarding non-injectivity of $\edgemap_n$.
We also have the following square
\[ \begin{tikzcd}
X_n \dar[two heads] \rar{\edgemap_n} & \prod X_1 \dar[two heads] \\
X_n^\flat \rar{\edgemap_n} & \prod X_1^\flat
\end{tikzcd} \]
but in general the bottom map will not be an injection.
Notice, though, that if $X_1 \to (X^\flat)_1$ is a bijection, then the bottom map will be an injection and we can conclude that $X^\flat$ is spiny.
In general we must iterate the construction.

\begin{definition}
We have a sequence of surjections
\[
	X = X^{(0)} \twoheadrightarrow X^{(1)} \twoheadrightarrow X^{(2)} \twoheadrightarrow \cdots  
\]
by defining $X^{(0)} \coloneqq X$ and $X^{(i+1)} \coloneqq (X^{(i)})^\flat$ for $i\geq 0$.
We let $X^{(\infty)}$ be the colimit of the sequence.
\end{definition}

If $U$ is a spiny symmetric set then 
\[
	\cdots \to \hom(X^{(i+1)}, U) \to \hom(X^{(i)}, U) \to \cdots \to \hom(X^{0},U) = \hom(X,U)
\]
is a sequence of bijections by \cref{lem flat}, so $X\to X^{(\infty)}$ induces a bijection
\[ 
\hom(X^{(\infty)},U) \cong \lim_{i\geq 0} \hom(X^{(i)}, U) \to \hom(X,U).\]

\begin{lemma}\label{lem infty spiny}
If $X$ is a symmetric set, then $X^{(\infty)}$ is a spiny symmetric set.
\end{lemma}
\begin{proof}
If $x\in X$, we write $[x]_k$ for the image of $x$ in $X^{(k)}$. 
Suppose $x,y \in X_n$ and suppose that $[x]_\infty, [y]_\infty \in X^{(\infty)}_n$ have the property that $\edgemap_n([x]_\infty) = \edgemap_n([y]_\infty) \in \prod X^{(\infty)}_1$.
For $1 \leq i \leq n$ we know $[\rho_{i-1,i}^*x]_\infty = [\rho_{i-1,i}^*y]_\infty$, so there exists a finite $k_i$ with $[\rho_{i-1,i}^*x]_{k_i} = [\rho_{i-1,i}^*y]_{k_i}$.
Take $k = \max k_i$, so that $\edgemap_n([x]_k) = \edgemap_n([y]_k)$.
But since $X^{(k+1)} = (X^{(k)})^\flat$, this implies that $[x]_{k+1} = [y]_{k+1}$, and hence $[x]_\infty = [y]_\infty$.
We conclude that $\edgemap_n \colon X^{(\infty)}_n \to \prod X^{(\infty)}_1$ is injective.
\end{proof}

\begin{example}
The following example shows that $X^{(m)}$ need not be spiny for any finite $m$, so iteration of $(-)^\flat$ need not stabilize in finitely many steps.
Consider the set of pairs $(a,f)$ where $a\in \{u,v\}$ and $f \colon [n] \to \mathbb{N} = \{0,1, \dots \}$ is a function landing in a three-element subset $\{ 0, k, k+1\}$ for some $k\in \mathbb{N}$.
Let $X_n$ denote the quotient under the identification $(u,f) \sim (v,f)$ whenever the image of $f$ is contained in a consecutive two-element set $\{k,k+1\}$.
The evident formula $\phi^*(a,f) \coloneqq (a, f\phi)$ turns $X$ into a symmetric set, which is pictured top left in \cref{fig symmetric set}. 
The symmetric set $X$ has the property that $X^{(m)}$ is not spiny for any finite $m$. 
One can show that $X^{(m)}$ is obtained by further identifying $(u,f) \sim (v,f)$ whenever $\max f(i) \leq m+1$, as indicated in the figure by the collapsing at each stage of the violet $2$-simplex and its companion underneath.
Passing to $m=\infty$, we have $(X^{(\infty)})_n$ is just the collection of $f$ as in the second sentence of this example; this is a symmetric subset of the nerve of a groupoid, hence is spiny.

\begin{figure}
\centering
\includegraphics[scale=0.94]{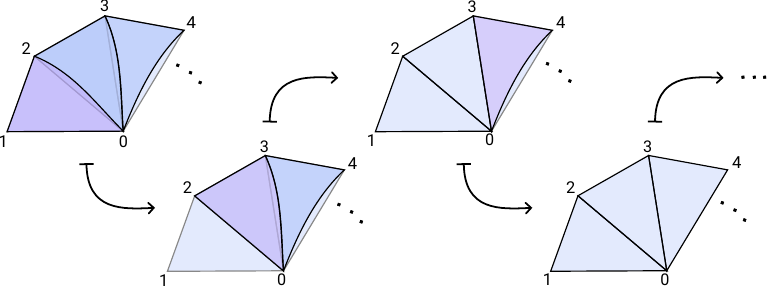}
\caption{Flattening a symmetric set}
\label{fig symmetric set}
\end{figure}
\end{example}

Since $X\to X^{(\infty)}$ induces a bijection $\hom(X^{(\infty)}, U) \cong \hom(X,U)$ for each spiny $U$, and the symmetric set $X^{(\infty)}$ is spiny, we obtain the following corollary. 

\begin{corollary}\label{cor spinification}
The category of (reduced) spiny symmetric sets is a reflective subcategory of (reduced) symmetric sets, via the reflection $X \to X^{(\infty)}$.  \qed
\end{corollary} 

\begin{remark}\label{rmk formal spinification}
We can recover \cref{cor spinification} from general principles without explicitly constructing $X^{(\infty)}$. 
It is straightforward to check directly that the full subcategory of spiny symmetric sets is closed in $\symset$ under products and subobjects, so \cite[16.9]{AHS:ACC} asserts\footnote{The category $\symset$ is a presheaf topos, so all of the formal requirements of this result hold. Such categories are wellpowered and complete, so strongly complete by \cite[12.5]{AHS:ACC}, and co-wellpowered by \cite[1.58]{AdamekRosicky:LPAC}. 
Finally, monomorphisms and extremal monomorphisms coincide in a topos, so `extremal subobjects' are the same as subobjects.} 
that spiny symmetric sets is an epireflective subcategory of $\symset$.
This means that it is a reflective subcategory and each reflection arrow is an epimorphism. 
The epimorphism property implies that this restricts to a reflection for reduced spiny symmetric sets inside reduced symmetric sets.
\end{remark}

\subsection{Consequences}\label{sec consequences}

Now that we have exhibited the categories of partial groups and partial groupoids as being equivalent to reflective subcategories of the presheaf category $\symset$, the following is standard.
This recovers results of \cite{Chermak:FL} and \cite{Salati:LCGRPG}.
\begin{theorem}\label{thm bicomplete}
The categories of partial groups and partial groupoids are complete and cocomplete. \qed
\end{theorem}
See, for instance, Proposition 3.3.9 and Proposition 4.5.15 of \cite{Riehl:CTC}.
Limits and colimits in $\symset$ exist and are constructed objectwise.
The limit of a diagram of (reduced) spiny symmetric sets in $\symset$ is again a (reduced) spiny symmetric set.
Colimits are more delicate: given a diagram $F \colon \mathcal{J} \to \pgpd$ of spiny symmetric sets, its colimit in $\symset$ is not typically spiny.
However, the reflection $(\colim_\mathcal{J} F)^{(\infty)} \in \pgpd$ will be a colimit in $\pgpd$ for the original diagram $F$.
If $F$ lands in partial groups instead, then $\colim_\mathcal{J} F \in \symset$ is often not reduced, in which case its reflection $(\colim_\mathcal{J} F)^{(\infty)} \in \pgpd$ will not be reduced either. 
Instead, $\mathcal{R}((\colim_\mathcal{J} F)^{(\infty)}) \cong (\mathcal{R}\colim_\mathcal{J} F)^{(\infty)}$ will be the colimit in $\pgrp$ of the original diagram of partial groups.

\begin{example}\label{ex colimits}
We collect several examples of colimits.
\begin{enumerate}[label=(\alph*), ref=\alph*]
\item Given a collection of partial groupoids, their coproduct as a symmetric set is spiny.
Hence their coproduct in the category of partial groupoids coincides with their coproduct in the category of symmetric sets.
\item Given a collection of partial groups, their coproduct as partial groups coincides with their coproduct in the category of reduced symmetric sets.
At each level, this is given by a wedge sum of pointed sets, so we will use $\vee$ for this coproduct.
\item If $\mathcal{J}$ is a connected category, then the colimit of any $\mathcal{J}$-shaped diagram of partial groups can be computed in either partial groups or partial groupoids with no change to the outcome. 
This includes, for instance, coequalizers, pushouts, and direct limits.\label{item connected colim}
\end{enumerate}
\end{example}

Theorem A.5 of \cite{Chermak:FL} states that certain colimits of partial groups exist and are preserved by the forgetful functor to pointed sets. 
While it is the case that the colimits exist, the following example shows that pushouts, for example, are not preserved by the forgetful functor.

\begin{example}\label{ex counterexample}
Let $A,B$ be two copies of the partial group $\fpg{2}$ from \cref{example reduced representable}, where in $B$ we rename $a_i$ to $b_i$, and let $T \cong \fpg{1} \vee \fpg{1}$ be the free partial group on two generators $t_1, t_2$. 
Consider the following pushout of partial groups.
\[ \begin{tikzcd}
T \rar{a_1,a_2} \dar[swap]{b_1,b_2} \ar[dr, phantom, "\ulcorner" very near end] & A \dar\\ 
B \rar & P
\end{tikzcd} \]
The elements $a_1$ and $b_1$ are identified in $P_1$, as are the elements $a_2$ and $b_2$.
The elements $\alpha = [a_1|a_2]$ and $\beta = [b_1|b_2]$ represent the same element in $P_2$, since $P_2 \to P_1 \times P_1$ is injective and $[a_1] = [b_1]$ and $[a_2] = [b_2]$ hold.
But then \[
[a_3] = [d_1\alpha] = d_1[\alpha] = d_1[\beta] = [d_1\beta] = [b_3]
\] in $P_1.$
On the other hand, we could consider the pushout of the underlying sets:
\[ \begin{tikzcd}
B_1 \dar{=}& T_1 \lar \rar \dar{=}& A_1\dar{=} \\
\{e, b_1, b_2, b_3, b_1^\dagger, b_2^\dagger, b_3^\dagger\}
& \{e, t_1, t_2, t_1^\dagger, t_2^\dagger \} \lar \rar & \{e, a_1, a_2, a_3, a_1^\dagger, a_2^\dagger, a_3^\dagger\}
\end{tikzcd} \]
where $a_3$ and $b_3$ are \emph{not} identified.
This pushout of sets has nine elements, while $P_1$ has seven.
Thus the forgetful functor to (pointed) sets does not preserve the pushout above, yet $B \leftarrow T \rightarrow A$ satisfies the hypotheses of Chermak's theorem.
\end{example}

In the example we just saw, the pushouts in $\pgrp$ and $\symset$ are different. 
This situation can be avoided when the partial group at the apex is actually a group.

\begin{example}\label{ex pushout along G}
Suppose we have injective partial group maps $G \hookrightarrow X$ and $G\hookrightarrow Y$ with $G$ a group.
Then in the pushout 
\[ \begin{tikzcd}
G \rar[hook,"F"] \dar[hook,"H"'] 
\ar[dr, phantom, "\ulcorner" very near end]
& X \dar[hook,"H'"] \\
Y \rar[hook,"F'"'] & Z
\end{tikzcd} \]
in $\symset$, we have that $Z$ is a partial group.
We already know that $Z$ is reduced by \cref{ex colimits}\eqref{item connected colim}, so we just need to check that it is spiny.
Since the underlying pushouts occur in $\set$, both $F'$ and $H'$ are injections, which implies that \[ (X_n \to Z_n \to \prod Z_1) = (X_n \hookrightarrow \prod X_1 \hookrightarrow \prod Z_1)\] is injective, and similarly for $Y_n \to \prod Z_1$.
Suppose we have $z,z'\in Z_n$ with $\edgemap_n(z) = \edgemap_n(z')$; we assume that $z = H'(x)$ and $z' = F'(y)$, otherwise we must already conclude that $z=z'$. 
Write $\edgemap_n(x) = [x_1 | \cdots | x_n]$ and $\edgemap_n(y) = [y_1 | \cdots | y_n]$.
Since $H'(x_i) = F'(y_i)$ and $F,H$ are injective, there exists a unique $g_i$ with $F(g_i) = x_i$ and $H(g_i) = y_i$.
Since $G$ is a group, $[g_1 | \cdots | g_n]$ is an $n$-simplex of $G$, mapping to $x$ and $y$.
Hence $z = H'(x) = F'(y) = z'$ are equal.
\end{example}

A \emph{filtered category} $\mathcal{J}$ is a non-empty category such that 1) if $a,b$ are two objects, then there are morphisms $a \to c \leftarrow b$ for some object $c$, and 2) if $f,g \colon a \to b$ are two parallel arrows, then there is some arrow $h\colon b \to c$ with $hf = hg$.
Alternatively, these are precisely those $\mathcal{J}$ such that for every finite category $\mathcal{I}$ and every functor $F \colon \mathcal{I} \times \mathcal{J} \to \set$, the canonical function
\[
  \colim_{j\in \mathcal{J}} \lim_{i\in \mathcal{I}} F(i,j) \to \lim_{i\in \mathcal{I}} \colim_{j\in \mathcal{J}} F(i,j)
\]
is a bijection of sets.
A filtered colimit is the colimit of a functor whose domain is a filtered category (so the above bijection expresses that filtered colimits commute with finite limits in $\set$ \cite[3.8.9]{Riehl:CTC}).
The standard example of such is a \emph{direct limit}, that is, a colimit of an ascending chain of morphisms 
\[
  A^0 \to A^1 \to A^2 \to \cdots.
\]
In algebraic settings (such as in the category of groups), filtered colimits are computed at the level of underlying sets.

\begin{proposition}\label{prop dir colim}
The category of partial groupoids is closed under filtered colimits in $\symset$.
\end{proposition}
Since filtered categories are connected, this implies (by \cref{ex colimits}\eqref{item connected colim}) that the category of partial groups 
$\pgrp$ is also closed in $\symset$ under filtered colimits.
\begin{proof}
This holds because filtered colimits commute with finite limits in $\set$.
In more detail, suppose $F \colon \mathcal{J} \to \pgpd$ is a filtered diagram.
Monomorphisms are characterized by a pullback condition as in the following square, and the $\mathcal{J}$-indexed diagram of pullbacks
\[ \begin{tikzcd}
F(j)_n \rar["\id"] \dar["\id"'] 
\ar[dr, phantom, "\lrcorner" very near start]
& F(j)_n  \dar["\edgemap_n"] \\
F(j)_n \rar["\edgemap_n"'] & \displaystyle \prod_{i=1}^n F(j)_1 
\end{tikzcd} \]
is sent to a corresponding pullback by $\colim_{\mathcal{J}}$. This implies that \[ \edgemap_n \colon \colim_{\mathcal{J}} F(j)_n \to  \colim_{\mathcal{J}} \prod_{i=1}^n F(j)_1 \cong \prod_{i=1}^n \colim_{\mathcal{J}} F(j)_1 \]
is a monomorphism.
\end{proof}

The category of partial groupoids is not closed under reflexive coequalizers, hence not closed under sifted colimits in $\symset$ (cf. \cite[4.2]{LackRosicky:NLT}).
Indeed, the pushout in \cref{ex counterexample} may be computed as the colimit of the following reflexive diagram,
\[ \begin{tikzcd}
A \vee T \vee B \rar[shift left=1.5] \rar[shift right=1.5] & A \vee B \lar
\end{tikzcd} \]
and we have already seen that the colimits in $\pgrp$ (or $\pgpd$) and $\symset$ differ.

\begin{corollary}
\label{cor locpres}
The categories of partial groups and partial groupoids are locally (finitely) presentable.
\end{corollary}
\begin{proof}
The category of spiny symmetric sets is closed in $\symset$ under directed co\-lim\-its\footnote{Directed colimits are a special case of filtered colimits, with the additional requirement that the indexing category $\mathcal{J}$ has at most one morphism between any two objects. 
A category has directed colimits if and only if it has filtered colimits, and for a functor out of such category, preservation of one type implies preservation of the other \cite[pp.13--16]{AdamekRosicky:LPAC}.} by \cref{prop dir colim}, as is the category of reduced spiny symmetric sets by invoking \cref{ex colimits}\eqref{item connected colim}.
The result is now a consequence of \cite[1.46]{AdamekRosicky:LPAC}.
\end{proof}

The following remark gives an elegant approach to the reflection theorem, as well as all other major results from this section, via \emph{orthogonality classes}.
Standard accounts are \cite[\S5.4]{Borceux:HCA1} and \cite[\S1.C]{AdamekRosicky:LPAC}.
We only sketch the argument.

\begin{remark}[Alternate proof of reflection theorem]\label{remark alt proof}
An object $x$ is \emph{orthogonal} to a map $f\colon a \to b$ if, for each map $g\colon a\to x$, there is a unique map $h\colon b \to x$ with $g = hf$. 
If $\Sigma$ is a collection of maps in a category $C$, let $C_\Sigma \subseteq C$ be the full subcategory spanned by those objects $x$ which are orthogonal to every map in $\Sigma$.
\begin{enumerate}[left = 0pt]
  \item Let $\spine{n} \subseteq \rep{n} = \hom(-,[n])$ consist of those $\phi\colon [m] \to [n]$ whose image lands in $\{i,i+1\}$ for some integer $i$. Write
  \[
    c_n \colon \rep{n} \amalg_{\spine{n}} \rep{n} \to \rep{n}
  \]
  for the collapse map from a pair of (symmetric) simplices glued along their spines.
  For $n\geq 2$ and $X \in \symset$, observe that $\edgemap_n \colon X_n \to \prod X_1$ is injective if and only if $X$ is orthogonal to $c_n$.
  \item A symmetric set is reduced if and only if it is orthogonal to $\varnothing \to \ast = \rep{0}$.
\end{enumerate}
Setting $\Sigma = \{ c_n \mid n \geq 2 \}$, we have $\symset_\Sigma \subset \symset$ is the category of spiny symmetric sets.
If $\Sigma' = \Sigma \cup \{ \varnothing \to \ast \}$ then $\symset_{\Sigma'} = \symset_\Sigma \cap \symset_{\{\varnothing \to \ast\}}$ is the category of reduced spiny symmetric sets.
\begin{enumerate}[left = 0pt,resume]
\item Deduce \cref{reflection theorem}, \cref{prop dir colim}, and \cref{cor locpres} from \cite[Theorem 1.39]{AdamekRosicky:LPAC}. \Cref{thm bicomplete} holds since locally presentable categories are bicomplete.
\end{enumerate}
\end{remark}

\bibliographystyle{amsalpha}
\bibliography{pgds}
\end{document}